\documentclass{elsart}

\usepackage{amsfonts}
\usepackage{amstext}
\usepackage{amssymb,amsmath}

\def\R{\mathbb{R}}

\def\P{\mathbb{P}}
\def\E{\mathbb{E}}
\def\WW{\mathbb{W}}

\def\0{\mathbf{0}}
\def\1{\mathbf{1}}
\def\cc{\mathbf{c}}
\def\dd{\mathbf{d}}

\def\r{\mathbf{r}}
\def\x{\mathbf{x}}
\def\X{\mathbf{X}}
\def\uu{\mathbf{u}}
\def\vv{\mathbf{v}}

\def\y{\mathbf{y}}
\def\z{\mathbf{z}}

\def\eps{\varepsilon}

\def\I{\mathbf{I}}

\def\M{\mathbf{M}}

\def\W{\mathbf{W}}
\def\A{\mathbf{A}}
\def\BB{\mathbf{B}}

\def\DD{\mathbf{D}}
\def\EE{\mathbf{E}}

\def\PPP{\mathbf{P}}

\def\tr{\mathtt{tr}\,}
\def\diag{\mathtt{diag}\,}
\def\Span{\mathtt{Span}\,}
\def\Var{{\mathtt{Var}}\,}
\def\Vol{{\mathtt{Vol}}}
\def\dist{{\mathtt{dist}}}

\def\part{\cal P}

\begin{document}

\begin{frontmatter}
\title {Modularity spectra, eigen-subspaces, \\and structure of weighted graphs}
\author{Marianna Bolla\corauthref{ca}}
\ead{marib@math.bme.hu}
\corauth[ca]{Research supported in part by the Hungarian 
National Research Grants OTKA 76481 and OTKA-KTIA 77778; further, by the 
T\'AMOP-4.2.2.C-11/1/KONV-2012-0001 project. Latter project has been supported
by the European Union, co-financed by the European Social Fund.}

\address{Institute of Mathematics,
Budapest University of Technology and Economics and \\
Inter-University Centre for Telecommunications and Informatics, Debrecen}

\begin{abstract} 
The role of the normalized modularity matrix in finding homogeneous 
cuts will be presented.
We also discuss the testability of the structural  
eigenvalues and that of the subspace spanned by the corresponding
eigenvectors  of this matrix.
In the presence of a spectral gap between the $k-1$ largest absolute value 
eigenvalues
 and the remainder of the spectrum, this in turn implies the testability of 
the sum of the inner variances of the $k$ clusters
that are obtained by applying the $k$-means algorithm for the appropriately
chosen vertex representatives.
\end{abstract}


\end{frontmatter}

\section{Introduction}
\label{intro}

The purpose of this paper is to summarize the spectral properties and
testability of the spectrum and spectral subspaces of the
normalized modularity matrix introduced in~\cite{Bol6} to find regular
vertex partitions. We will generalize the Laplacian based spectral
clustering methods to recover so-called volume regular cluster pairs such that 
the information flow between the pairs and within the clusters is 
as homogeneous as possible. For this purpose,
we take into consideration both ends of the normalized Laplacian spectrum, i.e.,
large absolute value, so-called structural eigenvalues of our normalized
modularity matrix introduced just for this convenience.

In Theorem~\ref{thk}, we estimate the constant of volume regularity in terms
of the gap between the structural and other eigenvalues, and the $k$-variance
of the optimal vertex representatives constructed by the eigenvectors
corresponding to the structural eigenvalues. Here we give a more detailed proof
of this statement than in~\cite{Bol7}.
This theorem implies that for a general edge-weighted graph,   
the existence of $k-1$ structural eigenvalues of the normalized modularity 
matrix, separated from 0, is
indication of a $k$-cluster structure such that the cluster-pairs are
volume regular with constant depending on the spectral gap and the above
$k$-variance. 
The clusters themselves can be recovered by applying
the $k$-means algorithm for the vertex representatives.
Hence, Theorem~\ref{thk} implies that spectral clustering of the vertices
into $k$ parts gives satisfactory partition in the sense of 
volume regularity.  

Furthermore, in Theorems~\ref{sajkonv} and~\ref{alterkonv}, we prove the
 testability of the structural eigenvalues and the corresponding
eigen-subspace of the normalized modularity matrix in the sense of~\cite{LovI}. 
In view of this, spectral clustering methods can be performed 
on a smaller part of the underlying graph and give good approximation 
for the cluster structure.

\section{Preliminaries}
\label{pre}

Throughout the paper, we use the general framework of an edge-weighted graph. 
Let $G=G_n =(V, \W)$ be an edge-weighted graph on vertex-set $V$ ($|V|=n$) and
$n\times n$ symmetric weight-matrix $\W$ of non-negative real entries and 
zero diagonal. 
We will call the numbers $d_i =\sum_{j=1}^n w_{ij}$ $(i=1,\dots ,n)$  
\textit{generalized degrees}, and the diagonal matrix  
$\DD = \diag (d_1 ,\dots ,d_n )$ \textit{degree matrix}.  
In this and the next section, without loss of generality, 
$\Vol (V) =1$ will be assumed, where
the \textit{volume} of the vertex-subset $U\subseteq V$ is  
$\Vol (U) =\sum_{i\in U} d_i$.
In the sequel, we only
consider connected graphs, which means that $\W$ is irreducible.

In~\cite{Bol6}, we defined the normalized version of the modularity matrix 
(introduced in~\cite{New}) as
$\M_D =  \DD^{-1/2} \W \DD^{-1/2} -\sqrt{\dd} \sqrt{\dd}^T$, 
where $\sqrt{\dd }=(\sqrt{d_1},\dots , \sqrt{d_n})^T$, 
and we called it \textit{normalized modularity matrix}.
The spectrum of this matrix is in the [-1,1] interval, and 0 is always an 
eigenvalue with unit-norm eigenvector $\sqrt {\dd}$. Indeed, in~\cite{Bol1}
we proved that  1 is a single eigenvalue of 
$\DD^{-1/2} \W \DD^{-1/2}$ with corresponding unit-norm eigenvector 
$\sqrt {\dd}$, provided our graph is connected.
This becomes a zero eigenvalue of $\M_D$ with the same 
eigenvector, whence 1 cannot be an eigenvalue of $\M_D$ if $G$ is connected.
In fact, the introduction of this matrix is rather technical, 
the spectral gap, further, Lemma~\ref{expmix} and 
Theorem~\ref{thk} can better be
formulated with it. It can also be obtained from the normalized Laplacian
by subtracting it from the identity and depriving of its trivial factor.
Normalized Laplacian was used for spectral
clustering in several papers (e.g.,~\cite{Azran,Bol1,Bol2,Chung1,Meila}), 
the idea of which  can be summarized by means of the 
spectral  decomposition of the 
normalized modularity matrix. We introduce the 
following notation: the \textit{weighted cut} between the vertex-subsets 
$X,Y\subseteq V$ is 
$w(X,Y) =\sum_{i\in X} \sum_{j\in Y} w_{ij}$.
We will frequently refer to the following facts.

\begin{itemize}
\item[{(a)}] 
The spectral decomposition of $\M_D$ solves
the following quadratic placement problem.
For a given positive integer $k$ ($1<k<n$), we want to minimize 
$ Q_k = \sum_{i<j} w_{ij} \| \r_i -\r_j \|^2 $
on the conditions
\begin{equation}\label{const}
 \sum_{i=1}^n d_i \r_i \r_i^T =\I_{k-1}  \quad \textrm{and} \quad
 \sum_{i=1}^n d_i \r_i =\0 
\end{equation}
where the vectors $\r_1 ,\dots ,\r_n$ are $(k-1)$-dimensional 
\textit{representatives} 
of the vertices, which form the row vectors of the $n\times (k-1)$ matrix $\X$.
Denote the eigenvalues of $\M_D$, in decreasing order, by 
$1>\lambda_1 \ge \dots \ge \lambda_n \ge -1$ with corresponding  
unit-norm, pairwise orthogonal eigenvectors $\uu_1 ,\dots ,\uu_{n}$.
In~\cite{Bol1}, we proved that the minimum of $Q_k$ subject to~(\ref{const}) 
is $k-1-\sum_{i=1}^{k-1} \lambda_i$ and
 is attained by the representation such that the optimum
vertex representatives $\r_1^* ,\dots ,\r_n^*$ are row vectors of the matrix 
$\X^* = (\DD^{-1/2} \uu_1 ,\dots ,\DD^{-1/2} \uu_{k-1} )$.
Instead of $\X$, the augmented $n\times k$ matrix $\tilde \X$ can as well be 
used, which is obtained from $\X$ by inserting the column $\x_0 =\1$ of all 1's.
In fact, $\x_0 =\DD^{-1/2} \uu_0$, where $\uu_0 =\sqrt{\dd }$ is the 
eigenvector corresponding to the eigenvalue 1 of $\DD^{-1/2} \W \DD^{-1/2}$.
Then
$$
Q_k = \tr (\DD^{1/2} {\tilde \X} )^T (\I_n - \DD^{-1/2} \W \DD^{-1/2} ) 
    (\DD^{1/2} {\tilde \X}) ,
$$
and minimizing $Q_k $ on the constraint~(\ref{const}) is equivalent to
minimizing the above expression subject to 
${\tilde \X}^T \DD {\tilde \X} =\I_k$.
This problem is the \textit{continuous relaxation} of minimizing 
$$
  Q_k (P_k )= \tr (\DD^{1/2} {\tilde \X} (P_k) )^T (\I_n - 
   \DD^{-1/2} \W \DD^{-1/2}  )(\DD^{1/2} {\tilde \X} (P_k ))
$$ 
over the set of
$k$-partitions $P_k =(V_1 ,\dots ,V_k )$ of the vertices such that $P_k$ is
planted into $\tilde \X$ in the way that the columns of ${\tilde \X} (P_k)$ 
are so-called \textit{normalized partition-vectors} belonging to $P_k$. Namely, 
the coordinates of the $i$th column are zeros,
except those indexing vertices of $V_i$, which are equal to 
$\frac1{\sqrt{\Vol (V_i )}}$ ($i=1,\dots ,k$).
In fact, this is the \textit{normalized cut} problem,
which is discussed in~\cite{Meila} for $k=2$,  further, 
in~\cite{Azran} and~\cite{Bol2} for a general $k$, 
and the solution is based on the above continuous relaxation.

\item[{(b)}] Now, let us maximize the
\textit{normalized Newman--Girvan modularity} of $G$ induced by $P_k$, 
defined in~\cite{Bol6} as
$$
 M_k ( P_k ) =  
 \sum_{a=1}^k \frac1{\Vol (V_a) }  \sum_{i,j\in V_a} (w_{ij} -d_i d_j ) 
  = \sum_{a=1}^k \frac{w(V_a ,V_a)}{\Vol (V_a )} -1  
$$
over the set ${\part}_k$ of the $k$-partitions of $V$. 
It is easy to see that 
$ M_k (P_k ) =k-1 -Q_k (P_k )$, and hence,
the above task has the same spectral relaxation as the normalized cut problem.
Let $ M_k =\max_{P_k \in {\part}_k} M_k (P_k ) $ denote the
maximum $k$-way normalized Newman-Girvan modularity of the weighted graph $G$.

\item[{(c)}] Finally, from the above considerations it is straightforward that
$M_k \le \sum_{i=1}^{k-1} \lambda_i$, or equivalently, the minimum normalized
$k$-way cut is at least the the sum of the $k-1$ smallest positive normalized
Laplacian eigenvalues. As for the minimum normalized $k$-way cut, 
in ~\cite{Bol2} we also gave an upper estimate  by constant times the
 sum of the $k-1$ smallest positive normalized
Laplacian eigenvalues, which constant depends 
on the so-called \textit{k-variance}
of the vertex representatives  defined in the following way.
\begin{equation}\label{kszoras}
 {S}_k^2 (\X ) =\min_{P_k \in {\part}_k }
  {S}_k^2 (\X , P_k ) =\min_{P_k =(V_1 ,\dots ,V_k )}
\sum_{a=1}^k \sum_{j\in V_a } d_j \| \r_j -{ \cc }_a \|^2 
\end{equation}
where ${\cc }_a =\frac1{\Vol (V_a ) } \sum_{j\in V_a } d_j \r_j $ is the
weighted center of cluster $V_a$  and $\r_1$, \dots ,$\r_n \in\R^{k-1}$ are 
rows of $\X$. (The augmented $\tilde \X$ would give the same $k$-variance.)
The constant of our estimation depended on
$S_k^2 ({\X}^* )$, and it was close to 1 if this $k$-variance of the optimum
$(k-1)$-dimensional vertex representatives was small enough.
Note that $ {S}_k^2 (\X , P_k )$ is the objective function of the weighted
$k$-means algorithm.
\end{itemize} 

In this way, we showed that large  positive eigenvalues of the normalized
modularity matrix are responsible for clusters with high
intra- and low inter-cluster densities. Likewise, maximizing $Q_k (P_k)$
instead of minimizing over ${\part}_k$,  small negative
eigenvalues of the normalized modularity
matrix are responsible for clusters with low
intra- and high inter-cluster densities (see~\cite{Bol6}). 
Our idea is that taking into
account eigenvalues from both ends of the normalized modularity
 spectrum, we can recover so-called regular cluster pairs. 
For this purpose, we use the notion of 
volume regularity to be introduced in the next section.

\section{Normalized modularity and volume regularity}
\label{normod}

With the normalized modularity matrix, the well-known Expander Mixing Lemma  
(for simple graphs see, e.g.,~\cite{Hoory}) is formulated 
for edge-weighted graphs  in the following way (see~\cite{Bol5}).
\begin{lem}\label{expmix} 
Provided $\Vol (V) =1$, for all $X,Y\subseteq V$,
$$
 | w (X, Y) - \Vol (X) \Vol (Y)| \le \| \M_D \|
 \cdot  \sqrt{\Vol (X) \Vol (Y)} ,
$$
where $\| \M_D \|$ denotes the spectral norm of the normalized modularity 
matrix of $G=(V, \W )$.
\end{lem}
Since the spectral gap of $G$ is  $1-\| \M_D \|$, a large spectral gap 
indicates small discrepancy as
a quasi-random property  discussed in~\cite{Chung}.
If there is a gap not at the ends of the spectrum,
we want to partition the vertices into clusters so that
a relation similar to the above property 
for the edge-densities between the cluster pairs 
would hold. For this purpose, we use a slightly modified version of the
volume regularity's notion introduced 
in~\cite{Alon}.
 \begin{defn} 
Let $G=(V, \W)$  be an edge-weighted graph with $\Vol (V) =1$.
The disjoint pair $A,B\subseteq V$ is
$\alpha$-volume regular
 if for all $X\subseteq A$, $Y\subseteq B$ 
we have
$$
| w (X, Y) -\rho (A,B) \Vol (X) \Vol (Y)| \le \alpha \sqrt{\Vol (A) \Vol (B)} ,
$$
where $\rho (A,B) =\frac{w(A,B)}{ \Vol (A) \Vol (B)}$ is the relative
inter-cluster density of $(A,B)$. 
\end{defn}
In the ideal $k$-cluster case, let us consider the following
generalized random simple graph model:
given the partition $(V_1 ,\dots ,V_k )$ of $V$ ($|V|=n$), 
vertices $i\in V_a$ and $j\in V_b$ are connected with
probability $p_{ab}$, independently of each other, $1\le a,b\le k$.
We can think of the probability $p_{ab}$ as the inter-cluster density
of the pair $(V_a ,V_b )$.
Since generalized random graphs can be viewed as edge-weighted graphs with a
special block-structure burdened with random noise, based on~\cite{Bol3},
we are able to give the following spectral characterization of them.
Fixing $k$, and tending with $n$ to infinity in
such a way that the cluster sizes grow at the same rate, 
 there exists a positive  number $\theta < 1$, independent
of $n$, such that for every $0<\tau <1/2$
there are exactly $k-1$  eigenvalues of 
$\M_D $ greater than $\theta -n^{-\tau}$,
while all the others are at most $n^{-\tau}$ in absolute value. Further, the
$k$-variance of the vertex representatives constructed by the $k-1$
transformed structural eigenvectors is ${\cal O}(n^{-2\tau})$,
and the cluster pairs are
$\alpha$-volume regular with any small $\alpha$, almost surely. Note that
generalized quasirandom graphs defined in~\cite{Lov} are deterministic
counterparts of generalized random graphs with the same spectral properties.

\begin{thm}\label{thk}
Let $G=(V, \W)$ be a connected edge-weighted graph on $n$ vertices, 
with generalized
degrees $d_1 ,\dots ,d_n$ and degree matrix $\DD$. Assume that 
$\Vol (V) =1$, and there are no dominant vertices, i.e.,
$d_i =\Theta (1/n )$, $i=1,\dots ,n$, as $n\to\infty$.   
Let the eigenvalues  of $\M_D$,
enumerated in decreasing absolute values, be 
$$
 1 \ge |\mu_1 | \ge \dots \ge |\mu_{k-1} | >\eps \ge |\mu_k | \ge
 \dots \ge |\mu_n |=0 . 
$$
The partition $(V_1,\dots ,V_k )$ of $V$ is defined so that it minimizes
the weighted k-variance $S_k^2 (\X^* )$ of the optimum vertex representatives
-- defined in~(\ref{kszoras}) --
obtained as row vectors of the $n\times (k-1)$ matrix $\X^*$ of column vectors
 $\DD^{-1/2} \uu_i$, 
where $\uu_i$ is the unit-norm eigenvector 
corresponding to $\mu_i$  $(i=1,\dots ,k-1 )$. 
Assume that there is a constant
$0<K\le \frac1{k}$ such that $|V_i |\ge Kn$, $i=1,\dots ,k$.
With the notation $s =\sqrt{S_k^2 (\X^* )}$, the $(V_i, V_j)$ pairs are
${\cal O} (\sqrt{2k} s +\eps)$-volume regular $(i\ne j)$ and for the
clusters $V_i$ $(i=1,\dots ,k )$  the following holds:
for all $X,Y\subset V_i$, 
$$
| w (X, Y) -\rho (V_i ) \Vol (X) \Vol (Y)| =
 {\cal O} (\sqrt{2k} s +\eps) \Vol (V_i ) ,
$$
where $\rho (V_i ) =\frac{w(V_i ,V_i )}{ \Vol^2 (V_i ) }$ is the relative
intra-cluster density of $V_i$. 
\end{thm}

Note that, in Section~\ref{pre}, we indexed the eigenvalues of $\M_D$
in non-increasing order and denoted them by $\lambda$'s. The set of all 
$\lambda_i$'s is the same as that of all $\mu_i$'s. Nonetheless, we need a
different notation for the eigenvalues indexed in decreasing order of 
their absolute values. Recall that 1 cannot be an eigenvalue of $\M_D$
if $G$ is connected. Consequently, $|\mu_1 |=1$ can be if and only if 
$\mu_1 =-1$, i.e., if $G$ is bipartite. For example, if the conditions of the 
above theorem hold with $k=2$ and $\mu_1 =-1$ ($|\mu_i |\le \eps$, $i\ge 2$), 
then our graph is a bipartite expander
discussed in~\cite{Alon0} in details.

For the proof we need the definition of the cut norm of a matrix 
(see e.g.,~\cite{Frieze}) and the 
relation between it and the spectral norm.
\begin{defn}\label{cutnorm}
The cut norm of the real matrix $\A$ with row-set $Row$ and 
column-set $Col$ is
$$
 \| \A \|_{\square } =\max_{R\subset Row, \, C\subset Col}
  \left| \sum_{i\in R} \sum_{j\in C} a_{ij} \right| .
$$
\end{defn}

\begin{lem}\label{cutspectral}
For the $m\times n$ real matrix $\A$,
$$
\| \A \|_{\square } \le \sqrt{mn} \| \A \|  ,
$$
where the right hand side contains the spectral norm, i.e. the largest singular
value of $\A$.
\end{lem}

\begin{pf}
$$
\begin{aligned}
 \| \A \|_{\square } &=\max_{\x \in \{ 0,1 \}^m,\, \y \in \{ 0,1 \}^n} 
 | \x^T \A \y | = 
\max_{\x \in \{ 0,1 \}^m ,\, \y \in \{ 0,1 \}^n} \left| 
(\frac{\x}{\| \x \|})^T \A 
(\frac{\y}{\| \y \|})\right| \cdot \| \x \| \cdot \|\y\| | \\
&\le \sqrt{mn}
\max_{\| \x \| =1, \, \|\y\|=1} | \x^T \A \y | =\sqrt{mn} \| \A \| ,
\end{aligned}
$$
since for $\x \in \{ 0,1 \}^m$, $\| \x \| \le\sqrt{m}$, and
for $\y \in \{ 0,1 \}^n$, $\| \y \| \le\sqrt{n}$. \qed
\end{pf}
The definition of the cut norm and the result of the above lemma naturally
extends to symmetric matrices with $m=n$, the spectral norm of which is the
absolute value of the maximum absolute value eigenvalue.

\begin{pf}(Theorem~\ref{thk}).
Recall that the spectrum of  $\DD^{-1/2} \W \DD^{-1/2}$ differs from that of
$\M_D$ only in the following: it contains the eigenvalue $\mu_0 =1$ with
corresponding unit-norm eigenvector $\uu_0 =\sqrt{\dd }$ instead of the 
eigenvalue 0 of $\M_D$ with the same eigenvector. 
If $G$ is connected, 1 is a single eigenvalue.
The optimum $(k-1)$-dimensional representatives of the vertices are row 
vectors   
of the matrix $\X^* =(\x_1^* ,\dots, \x_{k-1}^* )$, where 
$\x_i^* =\DD^{-1/2} \uu_i$
$(i=1,\dots ,k-1)$. The representatives can as well be regarded as 
$k$-dimensional
ones, as inserting the vector $\x_0^* =\DD^{-1/2} \uu_0 =\1$ will not change
the $k$-variance $s^2 =S_k^2 (\X^* )$.  
Assume that the minimum   $k$-variance is attained on the
$k$-partition $(V_1 ,\dots ,V_k )$ of the vertices. 
By an easy analysis of variance argument (see~\cite{Bol1}) it follows that
\begin{equation}\label{s}
 s^2 =\sum_{i=0}^{k-1} \dist^2 (\uu_i , F )  ,
\end{equation}
where $F =\Span \{ \DD^{1/2} \z_1 , \dots ,\DD^{1/2} \z_k \}$ 
with the so-called
normalized partition vectors $\z_1 ,\dots ,\z_k$ of coordinates
$z_{ji} = \frac1{\sqrt{\Vol (V_i )}}$ if $j\in V_i$ and 0, 
otherwise $(i=1,\dots ,k)$. 
Note that the vectors $\DD^{1/2} \z_1 , \dots ,\DD^{1/2} \z_k$ form an 
orthonormal system. By considerations proved in~\cite{Bol1}, 
we can find another orthonormal
system $\vv_0 ,\dots ,\vv_{k-1} \in F$ such that
\begin{equation}\label{s2}
  s^2 \le \sum_{i=0}^{k-1} \| \uu_i -\vv_i \|^2 \le 2s^2 
\end{equation} 
($\vv_0 =\uu_0$, since $\uu_0 \in F$).
We approximate 
the matrix $\DD^{-1/2} \W \DD^{-1/2} =\sum_{i=0}^{n-1} \mu_i \uu_i \uu_i^T$
by the rank $k$ matrix  $\sum_{i=0}^{k-1} \mu_i \vv_i \vv_i^T$ with the 
following accuracy (in spectral norm):
\begin{equation}\label{Frobenius}
 \left\| \sum_{i=0}^{n-1} \mu_i \uu_i \uu_i^T -
\sum_{i=0}^{k-1} \mu_i\vv_i\vv_i^T  \right\|
 \le \sum_{i=0}^{k-1} |\mu_i | \cdot \left\| \uu_i \uu_i^T -\vv_i \vv_i^T\right\|
 +\left\| \sum_{i=k}^{n-1} \mu_i \uu_i \uu_i^T \right\| 
\end{equation}
which can be 
estimated from above with 
$\sum_{i=0}^{k-1} \sin \alpha_i +\eps \le$
$\sum_{i=0}^{k-1} \| \uu_i -\vv_i \| +\eps \le \sqrt{2k} s +\eps$,
where $\alpha_i$ is the 
angle between $\uu_i$  and $\vv_i$, and for it,
$\sin \frac{\alpha_i}2 = \frac12 \| \uu_i -\vv_i \|$ holds, 
$i=0,\dots ,k-1$.

Based on these considerations and relation between the cut norm and the
spectral norm (see Lemma~\ref{cutspectral}),
the densities to be estimated in the defining formula  of
volume regularity can be written in terms of stepwise constant vectors in
the following way. The vectors $\y_i := \DD^{-1/2} \vv_i$ are stepwise 
constants on the partition $(V_1 ,\dots ,V_k )$, $i=0,\dots ,k-1$. The matrix 
$\sum_{i=0}^{k-1} \lambda_i \y_i \y_i^T $ is therefore a symmetric block-matrix 
on $k\times k$
blocks belonging to the above partition of the vertices. 
Let ${\hat w}_{ab}$
denote its entries in the $(a,b)$ block $(a,b=1,\dots ,k)$. 
Using~(\ref{Frobenius}), the rank $k$ approximation of the matrix $\W$ is
performed with the following  accuracy of the perturbation $\EE$:
$$
 \left\| \EE \right\| = \left\| \W - \DD (\sum_{i=0}^{k-1} 
 \mu_i \y_i \y_i^T ) \DD  \right\| =
 \left\| \DD^{1/2} ( \DD^{-1/2} \W \DD^{-1/2} -\sum_{i=0}^{k-1} \mu_i 
 \vv_i \vv_i^T ) \DD^{1/2} \right\| .
$$
Therefore, the entries of $\W$ -- for $i\in V_a$, $j\in V_b$ -- 
can be decomposed as
$
  w_{ij} = d_i d_j {\hat w}_{ab} +\eta_{ij} 
$,
where the cut norm  of 
the $n\times n$ symmetric error matrix $\EE =(\eta_{ij} )$  restricted to
$V_a \times V_b$ (otherwise it contains entries all zeros) 
and denoted by $\EE_{ab}$, is estimated as follows:
$$
\begin{aligned}
\| \EE_{ab} \|_{\square} &\le n\| \EE_{ab}\|\le n\cdot \| \DD_a^{1/2} \| \cdot 
(\sqrt{2k} s +\eps ) \cdot \| \DD_b^{1/2} \|  \\
 &\le n \cdot \sqrt{c_1 \frac{ \Vol (V_a )}{|V_a |} } \cdot
     \sqrt{c_1 \frac{ \Vol (V_b )}{|V_b |} } \cdot  (\sqrt{2k} s +\eps ) \\
 &=c_1 \cdot \sqrt{\frac{n}{|V_a|}} \cdot \sqrt{\frac{n}{|V_b|}} \cdot
 \sqrt{\Vol (V_a )} \sqrt{\Vol (V_b )}  (\sqrt{2k} s +\eps ) \\ 
 &\le
 c_1 \cdot \frac1{K} \sqrt{\Vol (V_a )} \sqrt{\Vol (V_b )} (\sqrt{2k} s +\eps )
 \\ &=c\sqrt{\Vol (V_a )} \sqrt{\Vol (V_b )}  (\sqrt{2k} s +\eps ) .
\end{aligned}
$$
Here the diagonal matrix $\DD_a$ contains the diagonal part of $\DD$ 
restricted to $V_a$, otherwise zeros, and 
the constant  $c$ does not depend on $n$.
Consequently, for $a,b=1,\dots ,k$ and $X\subseteq V_a$, $Y\subseteq V_b$:

\begin{eqnarray*}
 & &\left| w (X, Y) -\rho (V_a ,V_b ) \Vol (X) \Vol (Y)\right| =\\
 & & \left| \sum_{i\in X} \sum_{j \in Y} (d_i d_j {\hat w}_{ab}+
 \eta_{ij}) - 
\frac{\Vol (X) \Vol (Y )}{\Vol (V_a)  \Vol (V_b )}
\sum_{i\in V_a } \sum_{j \in V_b } (d_i d_j {\hat w}_{ab} +
\eta_{ij} )   \right| =\\
 & &\left| \sum_{i\in X} \sum_{j \in Y} \eta_{ij} - 
 \frac{\Vol (X) \Vol (Y )}{\Vol (V_a)  \Vol (V_b )}
 \sum_{i\in V_a } \sum_{j \in V_b } \eta_{ij} \right| \le 
 2c(\sqrt{2k} s +\eps )  \sqrt{\Vol (V_a ) \Vol (V_b )} ,
\end{eqnarray*}
that gives the required statement both in the $a\ne b$ and $a=b$ case.
\qed
\end{pf}


Note that
in the $k=2$ special case, due to a theorem proved in~\cite{Bol1}, the 
$2$-variance of the optimum 1-dimensional representatives can be directly
estimated from above by the gap between the two largest absolute value 
eigenvalues of $\M_D$,
and hence, the statement of Theorem~\ref{thk} simplifies, see~\cite{Bol5}.
For a general $k$, we can make the following considerations.

Assume that the normalized modularity spectrum 
(with decreasing absolute values) of $G =(V,\W )$  satisfies 
$$
 1\ge |\mu_{1}  |\ge \dots \ge |\mu_{k-1}| \ge \theta >\eps \ge 
|\mu_{k}| \ge \dots \ge |\mu_n |=0  .
$$
Our purpose is to estimate $s$ with the gap $\delta:= \theta -\eps$.
We will use the notation of the proof of Theorem~\ref{thk} and 
apply the results of~\cite{Bhat} for the perturbation of spectral 
subspaces of the symmetric matrices
$$
 \A =\sum_{i=0}^{n-1} \mu_i \uu_i \uu_i^T \quad \text{and} \quad 
 \BB =\sum_{i=0}^{k-1} \mu_i \vv_i \vv_i^T 
$$
in the following situation.
The subsets
$S_1 =\{ \mu_{k} ,\dots ,\mu_{n-1}\}$ and  
$S_2 =\{ \mu_{0} ,\dots ,\mu_{k-1} \}$ of the eigenvalues of 
$\DD^{-1/2} \W \DD^{-1/2}$
 are separated  by an annulus, where $\dist (S_1 ,S_2 ) =\delta >0$. 
Denote by $\PPP_A$ and $\PPP_B$ the projections onto the spectral subspaces
of $\A$ and $\BB$ spanned by the eigenvectors corresponding to the
eigenvalues in $S_1$ and $S_2$, respectively:  
$$
 \PPP_A (S_1 ) =\sum_{j=k}^{n-1} \uu_j \uu_j^T , \quad
 \PPP_B (S_2 ) =\sum_{i=0}^{k-1} \vv_i \vv_i^T .
$$ 
Then
Theorem VII.3.4 of~\cite{Bhat} implies that
\begin{equation}\label{bhatia}
 \| \PPP_A \PPP_B \|_F \le \frac1{\delta} \| \PPP_A (\A -\BB )\PPP_B  \|_F ,
\end{equation}
where $\| . \|_F$ denotes the Frobenius norm.
On the left hand side, 
$
  \| \PPP_A \PPP_B \|_F = \sqrt{\sum_{i=0}^{k-1}  \sin^2 \alpha_i } 
$, and in view of $\| \uu_i -\vv_i \| =2\sin \frac{\alpha_i}2$ and~(\ref{s2}), 
this is between $\frac{\sqrt{3}}{2} s$ and $s$. 
On the right hand side,
$$
 \PPP_A \A \PPP_B -\PPP_A \BB \PPP_B 
 = (\PPP_A \A ) \PPP_B -\PPP_A ( \PPP_B \BB )  
 = \sum_{i=0}^{k-1}  \sum_{j=k}^{n-1} (\mu_j -\mu_i )
 \uu_j^T (\uu_i -\vv_i ) 
 \uu_j \vv_i^T ,
$$
where the Frobenius norm of the rank 1 matrices $\uu_j \vv_i^T$ is 1, and the 
inner product  $\uu_j^T (\uu_i -\vv_i ) $ is the smaller if 
the $\uu_i$'s and the $\vv_i$'s are the closer $(i=1,\dots ,k-1)$.
Therefore, by the inequality~(\ref{bhatia}), 
$s$ is the smaller if $\delta $ is the larger and
the $|\mu_j -\mu_i |$ differences for $i=0,\dots ,k-1; \, 
j=k,\dots ,n-1$
are closer to $\delta$. If $|\mu_{k}| =\eps$ is small, then 
$|\mu_1 |,\dots,  |\mu_{k-1} |$ 
should be close to each other ($\mu_0 =1$ does not play an important
role because of $\uu_0 =\vv_0$).

\section{Testability of the normalized modularity spectrum and 
eigen-subspaces}
\label{test}

Authors of~\cite{LovI} defined the testability of simple graph parameters
and proved equivalent notions of this testability. 
They also anticipated 
that their results  remain valid if they consider weighted
graph sequences $(G_n )$  with edge-weights in the [0,1] interval and
 no dominant vertex-weights
$\alpha_i (G_n ) >0$ $(i=1,\dots ,n)$, i.e., 
$\max_i \frac{\alpha_i (G_n)}{\alpha_{G_n}} \to 0$ as $n\to\infty$,
where  $\alpha_{G_n}=\sum_{i=1}^n \alpha_i (G_n)$.
To this end, in~\cite{Bol4}, we slightly modified the
definition of a testable graph parameter for weighted graphs in the
following way. 
\begin{defn}\label{test}
A weighted graph parameter $f$ is testable if for every $\eps >0$
there is a positive integer $m<n$ such that if $G_n$ satisfies
$\max_i \frac{\alpha_i (G_n )}{\alpha_{G_n} } \le \frac1{m}$,
then 
$$
 \P (| f(G_n ) - f( \eta (m,G_n ) )| > \eps ) \le \eps ,
$$
where $\eta (m,G_n )$ is a random simple graph on  $m$ vertices
selected randomly from $G_n$  in the following manner:
 $m$ vertices of $G_n$ are selected  with replacement,
with respective probabilities proportional to the vertex-weights; 
given the selected vertex-subset, 
the edges come into existence conditionally independently, with
probabilities of the edge-weights. 
\end{defn}

By the above definition, a testable weighted graph parameter can be 
consistently estimated based on a fairly large sample. 
Based on the results of~\cite{LovI} for simple graphs, in~\cite{Bol4},
we established equivalent statements of this testability, from among which 
we will use the following.
\begin{fact}\label{facto}
Let $f$ be a testable weighted graph parameter. Then
for every convergent weighted graph sequence $(G_n)$, with
no dominant vertex-weights,
$f(G_n )$ is also convergent as $n\to\infty$.
\end{fact}

The notion of the convergence of a weighted graph sequence is defined 
in~\cite{LovI}, where the authors
also describe the limit object as a symmetric,
measurable function $W:[0,1]\times [0,1]\to [0,1]$, called \textit{graphon}.
The so-called \textit{cut distance} between the graphons $W$ and $U$ is
$
 \delta_{\square}  (W,U) =
\inf_{\nu}  \| W -U^{\nu} \|_{\square} 
$,
where the cut norm of the graphon $W$ is defined by
$$
 \| W  \|_{\square} =\sup_{S,T \subset [0,1]} 
  | \int_{S\times T}  W(x,y) \, dx \, dy  |  ,
$$
and the above infimum  is taken over all  measure preserving bijections
$\nu : [0,1] \to [0,1]$, while $U^{\nu}$ denotes the transformed $U$ after 
performing the same measure preserving bijection $\nu$ on both sides of the
unit square.
Graphons are considered modulo measure preserving maps,
and under graphon the whole equivalence class is understood. 
In this way, to a convergent weighted graph sequence $(G_n )$, there is
 a unique
limit graphon $W$  such that $\delta_{\square} (G_n ,W) \to 0$ as $n\to\infty$,
where  $\delta_{\square} (G_n ,W)$ is defined as $\delta_{\square} 
(W_{G_n} ,W)$ with the step-function graphon $W_{G_n}$ assigned to $G_n$
in the following way:
the sides of the unit square are divided into intervals
$I_{1} , \dots ,I_{n}$ of lengths $\alpha_1 (G_n ) /\alpha_{G_n} ,\dots ,
\alpha_n (G_n ) /\alpha_{G_n}$,
and over the rectangle $I_{i} \times I_{j}$ the stepfunction takes on the
value $w_{ij} (G_n )$. 

In~\cite{Bol4}, we proved the testability of some normalized and unnormalized
balanced multiway cut densities such that we imposed balancing conditions on
the cluster volumes. Under similar conditions, for fixed number of
 clusters $k$, the unnormalized and normalized multiway
cuts and modularities are also testable, provided our edge-weighted graph has
no dominant vertices. The proofs rely on statistical physics notions 
of~\cite{LovII}, utilizing the fact that the graph convergence implies the
convergence of the ground state energy (minimum of the energy function over
the set of $k$-partitions of vertices).
In~\cite{Re2}, the authors
showed that the Newman-Girvan modularity is an energy function (Hamiltonian),
and hence, testability of the maximum/minimum normalized modularities,
under appropriate balancing conditions, can be shown analogously. 
Here we rather discuss the testability of spectra and $k$-variances, because in
spectral clustering methods these provide us with polynomial time algorithms,
though only approximate solutions are obtained as analyzed in 
Section~\ref{pre}.

In Theorem 6.6 of~\cite{LovII}, the authors prove that the normalized spectrum 
of a convergent graph sequence also converges in the following sense.
Let $W$ be a graphon and $(G_n)$ be a sequence of weighted graphs with
uniformly bounded edge-weights tending to $W$. (For simplicity, we assume that
$|V (G_n )| =n$). Let $|\lambda_{n,1} |\ge |\lambda_{n,2} | \ge \dots \ge 
|\lambda_{n,n} |$
be the adjacency eigenvalues of $G_n$ indexed by their decreasing absolute 
values, and let $\mu_{n,i}=\lambda_{n,i}/n$ $(i=1,\dots ,n)$ be the normalized 
eigenvalues. 
Further, let $T_W$ be the $L^2 [0,1] \to L^2 [0,1]$ 
integral operator corresponding to $W$:
$$
 T_W f  (x) =\int_0^1 W(x,y) f(y) \, dy .
$$
It is well-known that his operator is self-adjoint and compact, and hence,
it has a discrete real spectrum, whose only possible point of accumulation is
the 0. Let $\mu_i (W )$ denote the $i$th largest absolute value eigenvalue
of $T_W$.  Then for every $i\ge 1$, $\mu_{n,i} \to \mu_i (W)$ as $n\to\infty$.
In fact, the authors prove a bit more (see Theorem 6.7 of~\cite{LovII}): 
if a sequence $W_n$ of uniformly bounded
graphons converges to a graphon $W$, then for every $i\ge 1$,
$\mu_{i} (W_n ) \to \mu_i (W)$ as $n\to\infty$.
Note that the spectrum of $W_G$ is the normalized spectrum of $G$,
together with countably infinitely many 0's. Therefore, the convergence
of the spectrum of $(G_n)$ is the consequence of that of $(W_{G_n})$.

We will prove that in the absence of dominant vertices, the normalized 
modularity spectrum is testable. 
To this end, both the modularity matrix and the graphon are related
to kernels of special integral operators, described herein.
Let $(\xi ,\xi' )$ be a pair of identically distributed real-valued 
random variables 
defined over the
product space ${\cal X}\times {\cal X}$ having a symmetric joint distribution 
$\WW$
with equal margins $\P$.
Assume that the dependence between $\xi$ and $\xi'$
is regular, i.e., their joint
distribution $\WW$ is absolutely continuous with respect to the
product measure $\P \times \P$, and let $w$ denote its Radon--Nikodym 
derivative, see~\cite{Renyi}. 
Let $H=L^2 (\xi )$ and $H' =L^2 (\xi' )$ be the Hilbert spaces of
random variables which are functions of $\xi$ and $\xi'$ and have zero 
expectation
and finite variance with respect to $\P$. Observe that 
$H$ and $H'$ are isomorphic Hilbert spaces with the covariance as inner 
product;
further, they are embedded as subspaces into the $L^2$-space defined similarly
over the product space. (Here $H$ and $H'$ are also isomorphic in the sense that
for any $\psi \in H$ there exists a $\psi' \in H'$ and vice versa, such that
$\psi$ and $\psi'$ are identically distributed.)

Consider the linear operator taking conditional expectation 
between $H'$ and $H$ with respect to the joint distribution.
It is an integral operator and will be denoted by $P_{\WW} : H' \to H$ as
it is a projection restricted to $H'$ and projects onto $H$.
To $\psi' \in H'$ the operator $P_{\WW }$ assigns $\psi \in H$ such that
$\psi = \E_{\WW } (\psi' \, | \, \xi )$, i.e.,
$$  
 \psi (x) = \int_{\cal Y} w(x,y) \psi' (y) \, \P (dy) , \quad x\in {\cal X}.
$$ 
If
$$
 \int_{\cal X} \int_{\cal X} w^2 (x,y) 
  \P(dx) \P (dy) <\infty ,
$$
then $P_{\WW }$ is a Hilbert--Schmidt operator, therefore compact and
has spectral decomposition
$$
 P_{\WW } =\sum_{i=1}^{\infty} \lambda_i \langle .,\psi'_i \rangle_{H'} 
 \psi_i  ,
$$
where for the eigenvalues $|\lambda_i |\le 1$ holds and the 
eigenvalue--eigenfunction equation looks like
$$
 P_{\WW } \psi'_i =\lambda_i \psi_i  \quad (i=1,2,\dots ),
$$
where $\psi_i$ and $\psi'_i$ are identically distributed, whereas their joint
distribution is $\WW$.
It is easy to see that $P_{\WW }$ is self-adjoint and it
takes the constantly 1 random variable of $H'$ into the constantly 1 random
variable of $H$; however, the $\psi_0 =1 ,\psi'_0=1$ pair 
is not regarded as a function pair with eigenvalue $\lambda_0 =1$, 
since they have no
zero expectation. More precisely, the kernel is reduced to $w(x,y)-1$.

\begin{thm}\label{sajkonv}
Let $G_n =(V_n , \W_n )$ be the general entry of a convergent sequence of 
connected edge-weighted 
graphs whose 
edge-weights are in [0,1] and 
the vertex-weights are the generalized degrees. Assume that
there are no dominant vertices. Let $W$ denote the limit graphon of the
sequence $(G_n )$, and let
$$
 1\ge|\mu_{n,1}  |\ge  |\mu_{n,2}| \ge \dots \ge |\mu_{n,n}| =0
$$
be the normalized modularity spectrum of $G_n$ (the eigenvalues are indexed
by their  decreasing absolute values).
Further, let $\mu_i (P_{\WW })$ is the $i$th largest absolute value eigenvalue 
of the  integral operator 
$P_{\WW } : L^2 (\xi' ) \to L^2 (\xi )$ 
taking conditional expectation with respect to
the joint measure $\WW$ embodied by the normalized limit graphon $W$,
and $\xi ,\xi'$ are identically distributed random variables
with the marginal distribution of their
symmetric joint distribution  $\WW$.
Then for every $i\ge 1$, 
$$
  \mu_{n,i} \to \mu_i (P_{\WW }) \quad \textrm{as} \quad  n\to\infty .
$$
\end{thm}

\begin{pf}
In case of a finite $\cal X$ (vertex set) 
we have a weighted graph, and we will show that the operator
taking conditional expectation with respect to the joint distribution
determined by the edge-weights
corresponds to its normalized modularity matrix.  

Indeed, let ${\cal X} =V$, $|V|=n$, and  $G_n =(V, \W )$ be an edge-weighted
graph on the $n\times n$ weight matrix of the edges $\W$ with entries 
$W_{ij}$'s; now, they do not necessarily sum up to 1.
(For the time being,
$n$ is kept fixed, so -- for the sake of simplicity -- 
we do not denote the dependence of $\W$ on $n$). 
Let the vertices be also weighted with special weights
$\alpha_i (G_n ) :=\sum_{j=1}^n W_{ij}$, $i=1,\dots ,n$. Then the
step-function graphon $W_{G_n}$ is such that  $W_{G_n}(x,y)=W_{ij}$ whenever
$x\in I_i$ and $y\in I_j$, where the (not necessarily contiguous) intervals 
$I_1 ,\dots ,I_n$ form a partition of [0,1] such that the length of $I_i$ is
$\alpha_i (G_n ) /\alpha_{G_n}$ $(i=1,\dots ,n)$.

Let us transform $\W$ into a symmetric joint distribution
$\WW_n$ over $V\times V$.
The entries $w_{ij} =W_{ij}/ \alpha_{G_n}$ $(i,j=1,\dots ,n)$ embody this
discrete joint distribution of random variables $\xi$ and $\xi'$ which are
identically distributed with marginal distribution $d_1 ,\dots ,d_n$,
where $d_i = \alpha_i (G_n ) /\alpha_{G_n}$ $(i=1,\dots ,n)$. With the
previous notation $H=L^2 (\xi )$,  $H'=L^2 (\xi' )$, the operator
$P_{\WW_n} : H'\to H$ taking conditional expectation is an integral operator
with now discrete kernel $K_{ij}=\frac{w_{ij}}{d_i d_j}$.
The fact that $\psi$, $\psi'$ is an eigenfunction pair of $P_{\WW_n}$ 
with eigenvalue $\lambda$ means that
\begin{equation}\label{eigen}
 \frac1{d_i} \sum_{j=1}^n w_{ij} \psi' (j) =
 \sum_{j=1}^n \frac{w_{ij}}{d_i d_j} \psi' (j) d_j  =\lambda \psi (i) ,
\end{equation}
where $\psi (j) =\psi' (j)$ denotes the value of $\psi$ or $\psi'$ taken  on
with probability $d_i$ (recall that $\psi$ and $\psi'$ are identically
distributed).
The above equation is equivalent to
$$
 \sum_{j=1}^n \frac{w_{ij}}{\sqrt{d_i} \sqrt{ d_j}} \sqrt{d_j} \psi (j) =
 \lambda \sqrt{d_i} \psi (i) ,
$$
therefore the vector of coordinates  $\sqrt{d_i} \psi (i)$
($i=1,\dots ,n$) is
a unit-norm eigenvector of the normalized modularity matrix with eigenvalue
$\lambda$ (note that the normalized modularity spectrum does not depend on the
scale of the edge-weights, it is the same whether we use  $W_{ij}$'s or
$w_{ij}$'s as edge-weights). 
Consequently, the eigenvalues of the conditional expectation 
operator are the same as the eigenvalues of the normalized modularity matrix,
and the possible values taken on by the eigenfunctions of the
conditional expectation operator are the same 
as the coordinates of the transformed eigenvectors of the 
normalized modularity matrix
forming the column vectors of the matrix  $\X^*$ of the optimal 
$(k-1)$-dimensional representatives, see Section~\ref{pre} (a).

Let $f$ be a stepwise constant
function on [0,1], taking on value $\psi (i)$ on $I_{i}$. 
Then $\Var \psi =1$ is equivalent to $\int_0^1 f^2(x) \, dx =1$.
Let $K_{G_n}$ be the stepwise constant graphon defined as $K_{G_n}(x,y)=K_{ij}$ 
for $x\in I_i$ and $y\in I_j$.
With this, the eigenvalue--eigenvector equation~(\ref{eigen}) looks like
$$
 \lambda f(x) =\int_0^1 K_{G_n} (x,y) f(y) \, dy .
$$

The spectrum of $K_{G_n}$ is the normalized modularity spectrum of $G_n$
together with countably infinitely many 0's 
(it is of finite rank, and therefore,
trivially compact), and
because of the convergence of the weighted graph sequence $G_n$, in
lack of dominant vertices, the sequence of graphons $K_{G_n}$ also converges.
Indeed, the $W_{G_n} \to W$ convergence in the cut metric means the convergence
of the induced discrete distributions $\WW_n$'s  to the continuous $\WW$.
Since  $K_{G_n}$ and  $K$ are so-called copula transformations of those
distributions, in lack of dominant vertices (this causes the convergence
of the margins) they also converge, which in turn implies the 
$K_{G_n} \to K$ convergence in the cut metric.

Let $K$ denote the limit graphon of  $K_{G_n}$ $(n\to\infty )$. This will be
the kernel of the integral operator taking conditional expectation with
respect to the joint distribution $\WW$. It is easy to see that this operator
is also a Hilbert--Schmidt operator, and therefore, compact.
With these considerations the remainder of the proof is analogous to the proof 
of  Theorem 6.7 of~\cite{LovII}, where the authors  prove
that if the sequence $(W_{G_n})$ of graphons converges  to the limit 
graphon $W$,
then both ends of the spectra of the integral operators, 
induced by $W_{G_n}$'s 
as kernels, converge to the ends of
the spectrum of the integral operator induced by $W$ as 
kernel. We apply this argument for the spectra of the integral operators
induced by the kernels $K_{G_n}$'s and $K$. \qed
\end{pf}

Note that in~\cite{LovSzeg3}, kernel operators are also discussed,
but not with our normalization.

\begin{rem}\label{corspek}
By Fact~\ref{facto}, provided there are no dominant
vertices,  Theorem~\ref{sajkonv} implies that
for any fixed positive integer $k$, the 
$(k-1)$-tuple of the largest absolute value eigenvalues of the normalized 
modularity matrix is testable.  
\end{rem}

\begin{thm}\label{alterkonv}
Assume that there are constants $0<\eps <\theta \le 1$ such that
the normalized modularity spectrum (with decreasing absolute values) of 
any $G_n$  satisfies 
$$
 1\ge|\mu_{n,1}  |\ge \dots \ge |\mu_{n,k-1}| \ge \theta >\eps \ge 
|\mu_{n,k}| \ge \dots \ge |\mu_{n,n} | =0 .
$$
With the notions of Theorem~\ref{sajkonv}, and assuming that there are
no dominant vertices of $G_n$'s, the subspace spanned by the 
transformed eigenvectors $\DD^{-1/2}\uu_1$, \dots ,$\DD^{-1/2}\uu_{k-1}$
 belonging to the $k-1$ largest absolute value
eigenvalues of the normalized modularity matrix of $G_n$ also converges to
the corresponding $(k-1)$-dimensional subspace of $P_{\WW}$. 
More precisely, if $\PPP_{n,k-1}$ denotes
the projection onto the subspace spanned by the transformed eigenvectors 
belonging to $k-1$ largest absolute value  eigenvalues  of the 
normalized modularity matrix of $G_n$, and $\PPP_{k-1}$ denotes the projection 
onto the  corresponding eigen-subspace of $P_{\WW}$, then 
$\|  \PPP_{n,k-1} - \PPP_{k-1} \| \to 0$ as $n\to\infty $
(in spectral norm).
\end{thm}

\begin{pf}
If we apply the convergence fact $\mu_{n,i} \to \mu_i (P_{\WW })$
for indices $i=k-1$ and $k$, we get that there will be a gap of order 
$\theta -\eps -o(1)$
between  $|\mu_{k-1} (P_{\WW} )|$ and $|\mu_k (P_{\WW} )|$ too.

Let $P_{\WW ,n}$ denote the $n$-rank approximation of $P_{\WW}$ (keeping its
$n$ largest absolute value eigenvalues, together with the corresponding
eigenfunctions) in spectral norm. The projection $\PPP_{k-1}$ ($k<n$)
operates on the eigen-subspace spanned by the eigenfunctions belonging to the
$k-1$ largest absolute value eigenvalues of $P_{\WW,n}$ in the same way as 
on the corresponding $(k-1)$-dimensional subspace determined by $P_{\WW}$. 
With these considerations, we
apply the perturbation theory of eigen-subspaces  with
the following unitary invariant norm: the trace- or Schatten-norm  of the
Hilbert--Schmidt operator $A$ is $\| A \|_{\tr } =(\sum_{i=1}^{\infty}
\lambda_i^4 (A ) )^{1/4}$. 
Our argument with the finite ($k-1$) rank
projections  is the following. Denoting by $P_{\WW_n}$ the integral
operator belonging to the normalized modularity matrix of $G_n$ 
(with kernel $K_{G_n}$ introduced in the proof of Theorem~\ref{sajkonv}),
$$
\begin{aligned}
 \| \PPP_{n,k-1} - \PPP_{k-1} \| &=\| \PPP_{n,k-1}^{\perp } \PPP_{k-1} \|
\le \| \PPP_{n,k-1}^{\perp } \PPP_{k-1} \|_{\tr} \\ 
 \le &\frac{c}{\theta -\eps -o(1)} \| P_{\WW_n} -P_{\WW,n} \|_{\tr} 
\end{aligned}
$$
with constant  $c$ that is at most $\pi /2$ (Theorem VII.3.2 of~\cite{Bhat}).
But
$$
  \| P_{\WW_n} -P_{\WW,n} \|_{\tr}  \le
 \| P_{\WW_n} -P_{\WW} \|_{\tr} + \| P_{\WW} -P_{\WW,n} \|_{\tr} ,
$$
where the last term tends to 0 as $n\to \infty$, since the tail of the
spectrum (taking the fourth power of the eigenvalues) of a Hilbert--Schmidt 
operator converges.
For the convergence of the first term we use Lemma 7.1 of~\cite{LovI}, which
states that the trace-norm of an integral operator can be estimated from 
above by four  times the cut norm of the corresponding kernel.
But the convergence in the cut distance of the corresponding kernels to zero
follows from the considerations made in the proof of Theorem~\ref{sajkonv}. 
This finishes the proof.  \qed
\end{pf}

\begin{rem}
As the $k$-variance depends continuously on the above subspaces
(see the expansion~(\ref{s}) of $s^2$ in the proof of Theorem~\ref{thk}), 
Theorem~\ref{alterkonv} implies the testability of the 
$k$-variance as well. 
\end{rem}

\section{Summary}

The above results suggest that 
in the absence of dominant vertices,
even the normalized modularity matrix of a  
smaller part of the underlying weighted graph, selected at random with an
appropriate procedure, is able to reveal its cluster structure.
Hence, the gain regarding the computational time of this spectral clustering 
algorithm is twofold: we only use a smaller part of the graph and the
spectral decomposition of its normalized modularity matrix runs in polynomial 
time in the reduced number of the vertices. 
Under the vertex- and cluster-balance conditions this method can give
quite good approximations for the multiway cuts and helps us to find the number
of clusters and identify the cluster structure. 
In addition, taking into account both the positive and negative, large
absolute value eigenvalues together with eigenvectors, regular cuts can also
be detected, as the investigated spectral characteristics give good
estimates for the volume regularity's constant of the cluster pairs by
Theorem~\ref{thk}.
Such regular cuts are looked for in social or biological networks, e.g.,
if we want to find equally functioning synapses of the brain.

\begin{ack}
We thank the anonymous referee for his/her constructive comments.
\end{ack}

\end{document}